\documentclass{amsart}
\usepackage{mathrsfs, amssymb,  amsmath, amsthm, enumerate}
\newtheorem{prop}{Proposition}
\newcommand{\fldr}{\ensuremath{\mathbb{R}}}
\newcommand{\sn}[1][n]{\ensuremath{\mathsf{S}_{#1}}} 
 
\newcommand{\gl}[1][n]{\ensuremath{\mathsf{GL}(#1)}} 
\newcommand{\trace}{\operatorname{tr }}
\newcommand{\GL}{\,\mathscr L\,} 
\newcommand{\GR}{\,\mathscr R\,} 
\newcommand{\GD}{\,\mathscr D\,} 
\newcommand{\GH}{\,\mathscr H\,} 
\newcommand{\GJ}{\,\mathscr J\,} 
\title{Quadrics {\em via}  Semigroups}
\author{V. N. Krishnachandran}
\address{%
Vidya Academy of Science and Technology\\ 
THRISSUR - 680 501, KERALA, INDIA.} 
\email{krishnachandran.v.n@vidyaacademy.ac.in}
\begin{document}

\begin{abstract}
Let $M_2$ be the semigroup of linear endomorphisms of a plane.
 We show that the 
space of idempotents in $M_2$ is a hyperboloid of one sheet, 
the set of semigroup-theoretic inverses of a nonzero singular
element in $M_2$  is a hyperbolic paraboloid, and the set of nilpotent
elements in $M_2$ is a right circular cone. 
\end{abstract}
\maketitle

This is the story of  the rediscovery of  classical 
three-dimensional  geometry, especially the geometry of
quadric surfaces,  while  studying the semigroup $M_2(\mathbb R)$ 
of linear endomorphisms of a real plane. 
One of the surfaces  that appears prominently in this context is the 
hyperboloid of one sheet, referred to as {\em spaghetti bundle} in  
\cite{Samu:88}. In this story the spaghetti presents itself as the set of 
idempotents in $M_2(\mathbb R)$, the cone emerges as the 
set of nilpotent elements and the hyperbolic paraboloid as 
the set of semigroup-theoretic  inverses of a singular element. 

This rediscovery was briefly announced in \cite{Kris:00}.  
Generalizations  of some of the ideas presented here to 
semigroups of linear endomorphisms of higher dimensional
vector spaces are discussed in \cite{Kris:03}.  
The little bit of semigroup theory quoted below is based 
on \cite{Clif+:61}.

\section{The Semigroup $M_2(\fldr)$}

A set $S$ together with an associative binary operation  in $S$  
is called a  {\em semigroup}.
An element $e$ in $S$ is called an  {\em idempotent}  if $e^2=e$.
If $X \subseteq S$, the set of idempotents in $X$ is
denoted by $E(X)$. In any semigroup $S$ we can define  
certain equivalence relations, 
denoted by $\GL$, $\GR$, $\GJ$, $\GD$ and $\GH$, 
and  called Green's relations.
Let  $S^1$ denote  $S$ if $S$ has an identity
element. Otherwise, let it denote $S$ with an identity element $1$
adjoined. For   $a,b \in S$ the first three are defined by 
\begin{alignat*}{2}
&a\GL b \quad &&\Leftrightarrow\quad aS^1=bS^1\\
&a\GR b \quad &&\Leftrightarrow\quad S^1a=S^1b\\
&a\GJ b \quad &&\Leftrightarrow\quad S^1aS^1=S^1bS^1
\end{alignat*}
and the remaining ones by
$\GD=\GL\circ\GR=\GR\circ\GL$,  $ \GH=\GL\cap \GR$.

If $a\in S$, the set of all elements in $S$ which are 
$\GL$-equivalent to $a$ is denoted
by $L_a$ and is  called the $\GL$-class containing $a$.
The notations $R_a, J_a, D_a, H_a$ have similar meanings.
These  are  the  Green classes in $S$.
Two elements $a,a'$ in a semigroup $S$ are called  {\em inverse elements}  if
$aa'a=a$, $a'aa'=a'$. A semigroup in which every element has an inverse 
is called a {\em regular semigroup}.  

A  well-known example of a regular semigroup is  $M_n(\mathbb K)$ 
(where $\mathbb K=\mathbb R$, or $\mathbb K=\mathbb C$)   of  linear
endomorphisms of an $n$-dimensional vector space $V$ over $\mathbb K$ 
under composition of mappings. 
The set \sn\ of singular elements in $M_n(\mathbb K)$ is  a 
regular subsemigroup of $M_n(\mathbb K)$. 
Treating functions as right operators  we see that, 
for $a,b \in M_n(\mathbb K)$,
$a\GL b$ if and only if  $a$ and $b$ have the same range, and,  
$a \GR b$ if and only if $a$ and $b$ have
the same null space. Further, $a \GD b $ if and only if $a$ 
and $b$ have the same rank, 
and,  $a\GJ b$ is equivalent to  $a\GD b$.

As already indicated, the semigroup having  special interest  to us is 
$M_2(\mathbb R)$, denoted by $M_2$ in the sequel.
 We shall represent elements of  $M_2$ 
as square matrices of order $2$ relative  to some fixed ordered orthonormal basis for $V$. 
Listing out  the entries in the elements of $M_2$ row-wise we get vectors in $\fldr^4$. 
In this way we may identify $M_2$ with $\fldr^4$. 
The usual inner product in $\fldr^4$ can be represented using the trace
function, defined by 
$\trace (x)=  x_1+x_4$. 
If $x=(x_1,\ldots,x_4)$ and $y=(y_1,\ldots, y_4)$ 
then 
$$
\langle x,y\rangle  = x_1y_1+\cdots+x_4y_4  = \trace(x^*y).
$$
where  $x^*$ is the transpose of $x$.
%
%
\section{Geometry of the Green Classes in $M_2$}
 $M_2$ has  three $\GD$-classes,  namely,
$D_0$,  $D_1$ and $D_2$, where $D_k$ is the set of endomorphisms 
of rank $k$. 
Obviously we have $D_0=\{0\}$ which is simply a point. From the well-known
fact  (p.168 \cite{Bori:85})   that 
the space $M(m,n,k)$ of $m\times n$ matrices of rank $k$ is a manifold of 
dimension $k(m+n-k)$ we  immediately deduce that $D_1$ is a three-dimensional submanifold 
of $\fldr^4$.   What this means is that sufficiently small neighborhoods of every point 
in $D_1$ `looks like' a three-dimensional euclidean space. 
Lastly $D_2$ is the set $\gl[2]$ of all invertible elements in $M_2$. It is well known that $\gl[2]$
is a  four-dimensional submanifold of $\fldr^4$. 

To describe the $\GL$- and $\GR$-classes in $M_2$, we require some
geometric terminology.
A line in a linear space $U$ is an affine subspace of $U$ generated
by two distinct points (that is, vectors) in $U$ and a plane $\mathcal P$ in $U$ is an affine
subspace   generated by three non-collinear points in $U$.
If  $\mathcal P$  passes through the origin in $U$ then the set  $\mathcal 
P\setminus \{0\}$  is  called a  {\em punctured plane}  in $U$. In a
similar way we may define  a {\em punctured line}  in $U$. 

If $0=a\in M_2$, then $L_a=R_a=\{0\}$. Also, if $a\in \gl[2]$ then $L_a=R_a=\gl[2]$. The
classes $L_a$ and $R_a$, when 
$a\ne 0$ and $a\not\in \gl[2]$, are   the nontrivial $\GL$- and
$\GR$-classes in $M_2$. 
If $e=\left[\begin{smallmatrix} 1& 0 \\ 0 & 0 \end{smallmatrix}\right]$, a simple argument involving
range and null space shows that
\begin{align*}
L_e &= \left\{   \left[ \begin{smallmatrix} \alpha & 0 \\ \beta & 0 \end{smallmatrix}\right]  \,:\,
\alpha, \beta\in \mathbb R, (\alpha,\beta)\ne (0,0)\right\},  \\
R_e &=\left\{   \left[ \begin{smallmatrix} \alpha & \beta \\  0  & 0 \end{smallmatrix}\right]  \,:\,
\alpha, \beta\in \mathbb R,(\alpha,\beta)\ne (0,0)\right\}. 
\end{align*}
These are obviously punctured planes. That this is  true for any 
$0\ne a\in \sn[2]$ can be easily verified.  
\begin{prop}\label{Prop.ii.2.2}
If $0\ne a\in \sn[2]$, then $L_a$ and $R_a$ are punctured planes lying in $\sn[2]$.
\end{prop}

Surprisingly, the converse of Proposition
\ref{Prop.ii.2.2} is also true, that is, every   punctured plane in $\sn[2]$ is a  nontrivial 
$\GL$-class or $\GR$-class in $\sn[2]$.
A nontrivial $\GH$-class, being  the intersection of a nontrivial $\GL$-class and a
nontrivial $\GR$-class,  is a punctured line in $\sn[2]$.

\section{Intersection of $\mathsf S_2$ with $x_1+x_4=\lambda$}

We have seen that $D_1$ is a three-dimensional manifold sitting in $\fldr^4$.
 To know more about  this manifold  we consider intersections of
$\sn[2]$ with hyperplanes in $M_2$. 
A hyperplane  in $M_2$ is the set  of all
$
x=\left[\begin{smallmatrix}x_1 &x_2\\
x_3&x_4\end{smallmatrix}\right]
$
satisfying an equation of the form
$$
\alpha x_1+\beta x_2+\gamma x_3+\delta x_4=\lambda
$$
where  $a=\left[\begin{smallmatrix}\alpha&\gamma\\
\beta&\delta\end{smallmatrix}\right] \in M_2$. This equation
can be expressed in the form
$
\trace(ax)=\lambda.
$
We denote this  hyperplane  by  $P(a;\lambda)$
and  its  intersection   with \sn[2]\   by $SP(a;\lambda)$.
We first consider  the special case   $SP(I; \lambda)$ 
where $I=\left[ \begin{smallmatrix} 1 & 0 \\ 0 & 1 \end{smallmatrix}\right]$.

\begin{prop}\label{Prop.ii.4.1}
If $\lambda\ne 0$ then $SP(I;\lambda)$ is a 
hyperboloid  of one sheet. Also,  $SP(I;0)$ is a right circular cone
with vertex at the origin in $M_2$. 
\end{prop}
\begin{proof}
Consider the following points in $P(I;\lambda)$:  
$$
O'(\tfrac\lambda 2,0,0\tfrac \lambda 2),\,\,
A(\tfrac \lambda 2 +\tfrac 1{\sqrt 2},0,0, \tfrac \lambda 2 - \tfrac 1 {\sqrt 2}), \,\,
B(\tfrac \lambda 2, \tfrac 1 {\sqrt 2}, \tfrac 1 {\sqrt 2}, \tfrac \lambda 2),\,\,
C(\tfrac\lambda 2, -\tfrac1{\sqrt 2}, \tfrac 1{\sqrt 2},\tfrac\lambda 2).
$$
Using the inner product and the induced norm in $\fldr^4$, we can see
that  the vectors
$\overrightarrow{O'A}$,  $\overrightarrow{O'B}$,  $\overrightarrow{O'C}$  
form an orthonormal system in $P(I;\lambda)$.
We choose $O'$ as the origin
and the directed lines joining $O'$ to the points $A,B,C$ as the coordinate
axes. We shall refer  to  the  coordinates of a point in $P(I;\lambda)$
relative to this  axes  as the Bell-coordinates  of the point.

Let the Bell-coordinates of any point  $Q(x_1,x_2,x_3,x_4)$ in $P(I;\lambda)$
be $(X,Y,Z)$ so that we have 
$$
\overrightarrow{O'P}
=X\overrightarrow{O'A}+Y\overrightarrow{O'B}+Z\overrightarrow{O'C}. 
$$
Since
$ \overrightarrow{O'P}  =\overrightarrow{OP}-\overrightarrow{OO'}$, 
etc.,    we must  have
\begin{equation}\label{Eq.ii.4.1}
x_1  = \tfrac \lambda 2+\tfrac X{\sqrt 2},\qquad x_2 = \tfrac{Y-Z}{\sqrt
2},\qquad 
x_3  = \tfrac{Y+Z}{\sqrt 2},\qquad\,  x_4 = \tfrac \lambda 2-\tfrac X{\sqrt 2}.
\end{equation}
If we substitute the 
above expressions  in the equation defining the space
\sn[2], namely,
$
x_1x_4-x_2x_3=0,
$
we get 
\begin{equation}\label{Eq.ii.4.2}  
\quad X^2+Y^2-Z^2=\tfrac{\lambda^2}{2}.
\end{equation}
If we consider any set of
numbers $(X,Y,Z)$ satisfying Eq.\eqref{Eq.ii.4.2}, then  using the relations in
Eq.\eqref{Eq.ii.4.1}, we can see that $(X,Y,Z)$ are the Bell-coordinates of a
point on $SP(I;\lambda)$.   Therefore
Eq.\eqref{Eq.ii.4.2} is the equation of $SP(I;\lambda)$
relative to the Bell-axes. When $\lambda\ne 0$, this equation  represents a
hyperboloid of one sheet (\S 64 \cite{Bell:60}).

When $\lambda=0$, Eq.\eqref{Eq.ii.4.2} reduces to 
$
X^2+Y^2-Z^2=0,
$
which   represents a right circular cone   with vertex at the origin and semi-vertical
angle $\tfrac {\pi}{4}$. Note that this cone lies in the hyperplane $P(I;0)$.
\end{proof}

By direct computation or otherwise one can easily see that
$SP(I;0)$ is the set of nilpotent elements in $M_2$.

We now  explore the relations between the geometrical properties
of  $SP(I;\lambda)$ and the algebraic properties of
$M_2$.
Obviously  the center  of $SP(I;\lambda)$ is $O'$ which  is a point on the line joining 
$O=\left[\begin{smallmatrix}0&0\\0&0\end{smallmatrix}\right]$ and
$I=\left[\begin{smallmatrix}1&0\\0&1\end{smallmatrix}\right]$. Since this line
is independent of the choice of $\lambda$, the centers  of the
hyperboloids $SP(I;\lambda)$, for various values of $\lambda$,
all lie on a fixed line.
The axis of rotation of the hyperboloid is given by $X=0$, $Y=0$. These equations
 produce a matrix of the form
\begin{equation*}
\begin{bmatrix}\frac \lambda 2&x_2\\-x_2&\frac \lambda 2\end{bmatrix}=
\begin{bmatrix}\frac \lambda  2&0\\0&\frac \lambda 2\end{bmatrix}+
\begin{bmatrix}0&x_2\\-x_2&0\end{bmatrix}.
\end{equation*}
The set of all such elements form a line through the center of the hyperboloid
parallel to the line formed by the set of skew-symmetric elements in
$M_2$.
The asymptotic  cone of $SP(I;\lambda)$ from its center
(\S 78\cite{Bell:60}) is $X^2+Y^2-Z^2=0$, which is  independent of $\lambda$.

From the above discussion we can construct an image of  $\sn[2]$.
If we imagine $\lambda$ as `time',   $\sn[2]$ can be thought of 
as a space generated by the `moving' hyperboloid
$SP(I;\lambda)$. Initially, that is when
$\lambda=0$,  we have a degenerate hyperboloid, namely, the right
circular cone $SP(I;0)$ with vertex at $O$ and axis  the line 
$\left\{\left[\begin{smallmatrix}0&\alpha\\
-\alpha&0\end{smallmatrix}\right]\,:\, \alpha\in \fldr\right\}$. As time
advances numerically, the center  of the hyperboloid moves along the line 
$\left\{\left[\begin{smallmatrix}\alpha&0\\
0&\alpha\end{smallmatrix}\right]\,:\, \alpha\in \fldr\right\}$.
The axis of rotation of the moving hyperboloid is a
line through the center  parallel to the axis of  the `initial' hyperboloid.
Further, at time $\lambda$, 
the radius of the principal circular section of the hyperboloid is
$\tfrac{|\lambda|} {\sqrt{2}}$. Clearly, this increases as $\lambda$
increases numerically. Thus, the hyperboloid $SP(I;\lambda)$ expands
simultaneously as it moves away from the origin.

\section{The Space of Idempotents}

Let  $E(k)$  be  the set of  idempotents  of rank $k$. Obviously,
$E(0)=\{0\}$ and $E(2)=\{I\}$. Also, it is easy to see  that $E(1)=SP(I;1)$.
 The next result follows from this. 

\begin{prop}\label{Th.ii.5.1}
The space $E(1)$ of idempotents of rank $1$ in $M_2$ is a hyperboloid of
one sheet.
\end{prop}

The generators   of $E(1)$ through $e$ 
 turn out to be the sets $E(L_e)$ and $E(R_e)$.
To see this it is enough to show that these sets are lines through $e$, for
if a line lies wholly on a conicoid it must belong to a system of
generating lines of the conicoid (\S\  97   \cite{Bell:60}).
That  $E(L_e)$ and $E(R_e)$  are lines follows  from 
$$
E(L_e)=e+(1-e)\sn[2]e \qquad E(R_e)=e+e\sn[2](1-e), 
$$
which can be  verified by direct computation.

Having found the generating lines on $E(1)$, we next determine how these 
are organized into two systems of generators  as  in
\cite{Bell:60}.  Let  $\mathbf L_1$ be    the family of all
lines of the form $E(L_e)$ and  $\mathbf L_2$,  the family of all lines of the form $E(R_e)$.   
Since $V$ is  $2$-dimensional, no two distinct members of $\mathbf L_1$ (or, of $\mathbf L_2$)
intersect. Also,  every member of $\mathbf L_1$ intersects every member
(except one) of $\mathbf L_2$, and vice versa. Thus  $\mathbf L_1$ and 
$\mathbf L_2$ are the two families of generating lines of $E(1)$.
Though it is not an intrinsic property of $E(1)$, 
it is interesting to see that  members of $\mathbf L_1$  are members of 
the $\lambda$-system  and those of $\mathbf L_2$ are members of 
the $\mu$-system of generators of  $E(1)$ (in the terminology of \cite{Bell:60}). 

The equation of the plane of the principal
circular section of $SP(I;1)$ is $Z=0$, which  is equivalent to 
$x_2=x_3$, which, in turn,   is equivalent to the statement that the element
$x\in M_2$ is symmetric.
Thus the principal circular section of $E(1)$ consists of the
symmetric  elements in $E(1)$.
\section{Intersections  of $\mathsf S_2$  \\ with Arbitrary Hyperplanes}
The first thing  we notice is that the
coordinates $(X,Y,Z)$ of any point  on $SP(a;\lambda)$ relative to 
some  rectangular Cartesian axes  in $P(a;\lambda)$ satisfy a second
degree equation. Hence $SP(a;\lambda)$ is a quadric surface in $P(a;\lambda)$.
The nature of this surface depends, naturally, on $a$ and $\lambda$. 
If $\lambda\ne 0$, then dividing the equation of the plane $P(a;\lambda)$ by
$\lambda$, we see that the planes  $P(a;\lambda)$ and $P(\tfrac
{1}{\lambda}a;1)$ coincide. Hence we  need consider only the special
surface $SP(a;1)$. 

Let $a$ be nonsinguar.  If  $x\in SP(a;1)$ then  $\trace(ax)=1$ and so $ax\in E(1)$.
Hence $x\in a^{-1}E(1)$. The converse is obvious and so $SP(a;1)=a^{-1}E(1)$.
Since $E(1)$ is a hyperboloid of one sheet and since the map $x\mapsto a^{-1}x$ 
is an affine map of $P(I;1)$ onto $P(a;1)$, $SP(a;1)$ must be a hyperboloid
of one sheet. 

\begin{prop}
If  $a$ is nonsingular $SP(a;1)$ is a hyperboloid of one sheet.
\end{prop}

Any regular semigroup $S$ is equipped with a natural partial order
 \cite{Namb:80} defined by 
$$
x\leq y \quad\text{if and only if}\quad R_x\leq R_y\text{ and }x=fy\text{ for some }f\in
E(R_x). 
$$
If  $a\in M_2$ is  nonsingular  then, it can be shown that 
 $SP(a;1)$ is the set of nonzero elements in $\sn[2]$
which are less than $a$ under the natural partial order in $M_2$. 

Let $a\in \sn[2]$.  If $x\in\sn[2]$ is an inverse of $a$ 
then we have $axa=a$ and so $ax\in E(1)$ which implies that $x\in SP(a;1)$.
The converse also holds. Thus, if $0\ne a\in \sn[2]$, then $SP(a;1)$ is the
set of inverses of $a$.   
To understand the geometry of  $SP(a;1)$, we consider the special case 
$
e = \left[\begin{smallmatrix}1 & 0\\
0 & 0 \end{smallmatrix}\right]
$.
Clearly, $x\in P(e;1)$ if and only if $x_1=1$. 
Let $O'$,
$A$, $B$, $C$ denote $(1,0,0,0)$, $(1,1,0,0)$, $(1,0,1,0)$ and $(1,0,0,1)$ 
respectively. Then the vectors 
$\overrightarrow{O'A}$,  $\overrightarrow{O'B}$,  $\overrightarrow{O'C}$  
form an orthonormal basis in $P(e;1)$. If  $Q$ denotes $(1,x_2,x_3,x_4)$ 
then
$$
\overrightarrow{O'Q}=x_2 \overrightarrow{O'A} +  
x_3  \overrightarrow{O'B} + x_4 \overrightarrow{O'C}
$$
so that the coordinates $(X,Y,Z)$ of $Q$ relative to this system of axes are 
$X = x_2$, $Y=x_3$, $Z=x_4$.  Now $x$  is in $SP(e;1)$ if and only if $x_1=1$  and $x_4=x_2x_3$ and so
$SP(e;1)$ is determined by the cartesian equation 
$
XY=Z
$
which  is the equation of a hyperbolic paraboloid. This argument can be modified for arbitrary
$a\in \sn[2]$ to yield the following result.

\begin{prop}\label{Prop.ii.6.9}
Let $a$ be singular. Then $SP(a;1)$ is the set of inverses of $a$ and is a hyperbolic
paraboloid.
\end{prop}

We next consider intersections of $\sn[2]$ with 
hyperplanes passing through the origin. These are sets of the form 
$SP(a;0)$. If  $a$ is  nonsingular then  $SP(a;0)=a^{-1}SP(I;0)$.
Since $SP(I;0)$ is a cone with vertex at the origin,  $SP(a;0)$ is also a cone.
If  $a$  is  singular and nonzero then using some technical arguments, we can 
 show   that  $SP(a;0)\setminus\{0\}$ is the union
of an  $\GL$-class and an  $\GR$-class in \sn[2], that is, a union of two punctured planes.

\end{document}